\documentclass[12pt]{amsart}
\usepackage{latexsym,amssymb,amsmath}

\newcommand{\R}{{\mathbb R}}
\newcommand{\T}{{\mathbb T}}
\newcommand{\Z}{{\mathbb Z}}

\newcommand{\PP}{{\mathcal {P}}}
\newcommand{\CC}{{\mathcal {C}}}
\newcommand{\DD}{{\mathcal {D}}}

\newtheorem{theorem}{Theorem}[section]
\newtheorem{lemma}[theorem]{Lemma}
\newtheorem{proposition}[theorem]{Proposition}

\begin{document}

\title[  relations between bilinear multipliers on
$  \R^n, \mathbb{T}^n$ and $\mathbb{Z}^n$] {relations between bilinear
multipliers on $ \R^n, \mathbb{T}^n$ and $\mathbb{Z}^n$}
\author[Bose]{Debashish Bose}
\address{Independent researcher, Kanpur, India}
\email{debashishb@wientech.com}

\author[Madan]{Shobha Madan}
\address{Department of Mathematics and
 Statistics, I.I.T. Kanpur, India}
\email{madan@iitk.ac.in, parasar@iitk.ac.in, saurabhk@iitk.ac.in}

\author[Mohanty]{Parasar Mohanty}

\author[Shrivastava]{Saurabh Shrivastava}


\begin{abstract}
In this paper we prove the bilinear analogue of de Leeuw's result for
periodic bilinear multipliers and some Jodeit type extension results for
bilinear multipliers.

\end{abstract}

\maketitle

\section{\bf Introduction}
 In their much acclaimed work Lacey and Thiele \cite{lt1}, \cite{lt2} proved the boundedness of the bilinear Hilbert
 transform and established a long standing conjecture of A.P. Calder\'on. Since then the study of bilinear multiplier
 operators which commute with simultaneous translations have attracted a great deal of attention. For a
 comprehensive survey we would like to refer the interested
  reader to the article of Grafakos and Torres \cite{gtor}.
\medskip

One of the important themes  of study of $L^p$ multipliers is about
the relationship between multipliers on the classical Euclidean groups
$\R^n,\;\T^n,\;\Z^n$. de Leeuw \cite{dlw} studied the restrictions
of $L^p$ multipliers on $\R^n$ to $\T^n$. These kind of relations
 between  bilinear multiplier operators
defined on $\R$ and $\T$ have appeared in the work of  Fan and Sato
\cite{fs}, Blasco, Carro and Gillespie \cite{bcg} and Grafakos
\cite{ghoz}.

  Also, in the same paper de Leeuw observed that periodic multipliers on
$\R$ are precisely the ones which are multipliers on $\Z$ and vice
versa. In section 3 we investigate the bilinear analogue of these
results.

   The extension question from $\T^n$ to $\R^n$ was not fully
   explored in de Leeuw's paper. However, in \cite{jod} Jodeit
   addressed some natural extensions. A function of $\Z^n$ is extended to
   $\R^n$ by forming the sum of integer translates of a suitable
   function
   $\Lambda$. i.e.$$\Psi(\xi)=\sum\limits_{k\in\Z^n}\phi(k)\Lambda
   (\xi-k).$$
   (For $n=1$ if  $\Lambda$ takes the value $1$ at zero and has support in $[0,1)$
   then $\Psi$ and $\phi$ agree at integers). In \cite{sp} Madan and Mohanty addressed the bilinear analogue of
this using transference techniques. They have shown that the
piecewise constant extension of a bilinear multiplier symbol
$\phi(n,m)$ of an operator
 $\PP_\phi:L^{p_1}(\T) \times L^{p_2}(\T) \rightarrow L^{p_3}(\T)$,  where $\frac{1}{p_1}+ \frac{1}{p_2}=
 \frac{1}{p_3}<1$, gives a bilinear multiplier on $\R$. In section \ref{jext}, we give other examples of $\Lambda$
  as in Jodeit which complements the results of \cite{sp}. %and consequently extend the above mentioned result in \cite{sp}.
Further the results hold for the entire admissible range of exponents $p_1, p_2$ and $p_3$.

\medskip

In Section \ref{prelim}, we give
basic definitions and notation.

\section{\bf Preliminaries}\label{prelim}
\medskip

Let $\mathcal{S}(\R)$ be the space of Schwartz class functions with
the usual topology and let $\mathcal{S'}(\R)$ be
 its dual space. We say
that a triplet $(p_1, p_2, p_3)$ is H\"{o}lder related if
$\frac{1}{p_1}+ \frac{1}{p_2}= \frac{1}{p_3}$ ,where $p_1,p_2 \geq
1$ and $p_3 \geq \frac{1}{2}$.
\medskip

For $f, g \in \mathcal{S}(\R)$ the  Bilinear Hilbert transform   is given by
$$ H(f,g)(x):= p.v. \int_\R f(x-t) g(x+t) \frac{dt}{t} $$
and  has the following alternative expression  : $$H(f, g) (x)=
\int_{\R}\int_{\R}\hat{f}(\xi) \hat{g}(\eta)(-i )sgn(\xi - \eta)e^{2
\pi i x(\xi+\eta)} d\xi d\eta . $$

In \cite{lt1}, \cite{lt2}, Lacey and Thiele  proved the boundedness of the above
operator $H$ from $L^{p_1}(\R) \times L^{p_2}(\R) \rightarrow
L^{p_3}(\R)$ for the H\"{o}lder related triplet $(p_1, p_2,
p_3)$, where $1 < p_1, p_2 \leq \infty $ and $\frac {2}{3}< p_3
<\infty  $.

\medskip

It is known that for any continuous bilinear operator $\CC : \mathcal{S}(\R) \times \mathcal{S}(\R) \rightarrow \mathcal{S'}(\R)$,
which commutes with simultaneous translations there exists a symbol $\psi_{\CC}(\xi,\eta)$ such that for $f, g \in \mathcal{S}(\R)$

$$\CC(f, g) (x)= \int_{\R}\int_{\R}\hat{f}(\xi) \hat{g}(\eta)\psi_{\CC}(\xi,\eta)e^{2 \pi i x(\xi+\eta)} d\xi d\eta $$

In the distributional sense we can write
$$\CC(f, g) (x)= \int_{\R}\int_{\R} f(x - u) g(x-v) K_{\CC}(u, v) du dv $$
where $\hat{K_{\CC}}= \psi_{\CC}$ (in the sense of distributions).
\medskip

Unlike in the linear case, the boundedness of the symbol $\psi_{\CC}$ is not known. In this article we will be
dealing with bounded symbols only.
\medskip
For $\psi \in L^{\infty}(\R^2)$ and $f, g \in \mathcal{S}(\R)$, we write
\begin {eqnarray}\label{symbol side}
\CC_{\psi}(f, g) (x)= \int_{\R}\int_{\R}\hat{f}(\xi) \hat{g}(\eta)\psi(\xi,\eta)e^{2 \pi i x(\xi+\eta)} d\xi d\eta
\end {eqnarray}\\
If for all $ f,g \in \mathcal{S}(\R)$ the bilinear operator $\CC_{\psi}$ satisfies
\begin {eqnarray}\label{bilinear}||\CC_{\psi}(f,g)||_{p_3} \leq c ||f||_{p_1} ||g||_{p_2},\end {eqnarray}
where $c$ is a constant independent of $f$, $g$ , then  we say that
$\CC_{\psi}$ is a bilinear multiplier operator associated with the
symbol $\psi$ for the triplet $(p_1, p_2, p_3)$. The set of all
bounded bilinear multipliers for the triplet $(p_1, p_2, p_3)$ will
be denoted by $M_{p_1,p_2}^{p_3}(\R)$. For $p_3 \geq 1$,
$M_{p_1,p_2}^{p_3}(\R)$ becomes a Banach space under the operator
norm, whereas for $p_3 < 1$ it forms a quasi Banach space. We will use the notation $\|.\|$ for the operator norm and
for convenience we will not attach any $p$ with it. It will be understood from the context.

\medskip

Bilinear multiplier operators on $\T$ and $\Z$ can be defined similarly. We say that the operator $\PP$ defined by
$$\PP(F, G)(x):=\Sigma_n \Sigma_m \hat{F}(n)\hat{G}(m)\phi(n,m) e^{2\pi i x(n+m)}; $$
is bounded from $L^{p_1}(\T)\times L^{p_2}(\T) \rightarrow$
$L^{p_3}(\T)$ if for some constant $c$ and for all trigonometric
polynomials $F, G $  we have
$$||\PP(F,G)||_{p_3} \leq c ||F||_{p_1} ||G||_{p_2}$$

Similarly on $\Z$ the operator
$$\DD(a, b)(l)= \int_{\T}\int_{\T}\hat{a}(\theta) \hat{b}(\rho)\psi(\theta,\rho)e^{2 \pi i l(\theta+\rho)} d\theta d\rho$$
is said to be bounded from $l^{p_1}(\Z)\times l^{p_2}(\Z) \rightarrow l^{p_3}(\Z)$, if for some constant $c$ and for all
finite sequences $a, b$ we have
$$||\DD(a,b)||_{l^{p_3}} \leq c ||a||_{l^{p_1}} ||b|_{l^{p_2}}$$

The space of bounded bilinear multipliers on $\T$ and $\Z$ for the
triplet $(p_1, p_2, p_3)$ will be denoted by $M_{p_1,p_2}^{p_3}(\T)$
and $M_{p_1,p_2}^{p_3}(\Z)$ respectively.

\section {\bf Periodic bilinear multipliers}\label{periodic}
\medskip
%\subsection {\bf Proof of Theorem (\ref{periodic symbols})}

Let $ \psi({\xi,\eta}) \in M_{{p_1},{p_2}}^{p_3}(\R)$ be a periodic
function with period one in both variables, i.e.
$\psi({\xi,\eta})=\psi({\xi+1,\eta})=\psi({\xi,\eta+1}).$ A natural
question that arises is whether $\psi({\xi,\eta}) \in
M_{{p_1},{p_2}}^{p_3}(\Z)$. In \cite{blasco} Blasco proved a
partial result in this direction. Conversely, given $\psi \in
M_{{p_1},{p_2}}^{p_3}(\Z)$ one can ask whether $ \psi({\xi,\eta})
\in M_{{p_1},{p_2}}^{p_3}(\R)$. For linear multipliers see \cite{auscher}. We address these questions here for the entire admissible range of exponents. In particular, we show that

\begin{theorem}\label{periodic symbols} Let $\psi \in \, L^\infty(\R^2)$ be
a $1$-periodic function in both variables. Then $\psi \in
M_{{p_1},{p_2}}^{p_3}(\R)$ if and only if  $\psi \in
M_{{p_1},{p_2}}^{p_3}(\Z)$, where the triplet $(p_1, p_2, p_3)$ is
H\"{o}lder related.\end{theorem}
\noindent
{\bf Proof}: First we will prove that if $\psi\in M_{{p_1},{p_2}}^{p_3}(\R)$ then
$\psi\in M_{{p_1},{p_2}}^{p_3}(\Z)$. Let $\Phi\in \mathcal{S}(\R)$
be such that $supp(\Phi)\subseteq[0,1], \, 0\leq\Phi(x)\leq1$ and
$\Phi(x)=1, \forall \, x\in[\frac{1}{4},\frac{3}{4}]$. If $\{a_k\}$
and $\{b_l\}$ are two  sequences with finitely many non-zero terms,
we define two functions $f_a,g_b$ as follows
$$f_a(x):= \sum\limits_k a_k \Phi(x-k)$$ and $$g_b(x):= \sum\limits_l b_l \Phi(x-l).$$
It is easy to see that $||f_a||_{L_{p_1}({\R})} \leq ||a||_{l_{p_1}({\Z})}$
and $||g_b||_{L_{p_2}({\R})} \leq ||b||_{l_{p_2}({\Z})}$. Then

$$\CC_\psi(f_a,g_b)(x)= \int_{\R}\int_{\R} \psi(\xi,\eta)  \hat{f_a}(\xi)
\hat{g_b}(\eta) e^{2\pi i (\xi+\eta)x}  d\xi d\eta$$

$$= \int_0^1\int_0^1 \psi(\xi,\eta)e^{2\pi i (\xi+\eta)x} \sum\limits_m \sum\limits_n\hat{f_a}(\xi+m)\hat{g_b}(\eta+n)
e^{2\pi i (m+n)x} d\xi d\eta $$

$$= \int_0^1\int_0^1 \psi(\xi,\eta)e^{2\pi i (\xi+\eta)x} \sum\limits_m\hat{f_a}(\xi+m)e^{2\pi i mx}
\sum\limits_n\hat{g_b}(\eta+n)e^{2\pi i nx} d\xi d\eta $$

%$$=\int_0^1\int_0^1 \psi(\xi,\eta)e^{2\pi i (\xi+\eta)x}\left(\sum_m f_a(x+m)e^{-2\pi i \xi(x+m)}\right)\left(\sum_n g_b(x+n)e^{-2\pi i \eta(x+n)}\right) d\xi  d\eta$$

$$= \int_0^1\int_0^1 \psi(\xi,\eta)\left(\sum\limits_m f_a(x+m)e^{-2\pi i \xi m}\right)\left(\sum\limits_n g_b(x+n) e^{-2\pi i \eta  n}\right) d\xi  d\eta $$

where we have used the Poisson summation formula in the last step.
\medskip

For $x\in[j+\frac{1}{4},j+\frac{3}{4}] = I_j$, we can write $x = [x]+
(x')$, where $(x')$ is the fractional part of $x$.Then
\begin{eqnarray*}
\sum\limits_m f_a(x+m)e^{-2\pi i \xi m}
&=& \sum\limits_m\sum\limits_k a_k \Phi(j+(x')+m-k)e^{-2\pi i \xi m}\\
&=& \sum\limits_m\ a_{j+m}e^{-2\pi i \xi m}\\
 &=&\sum\limits_m\ a_m e^{-2\pi i \xi (m-j)}\\
 &=&  \hat{a}({\xi}) e^{2\pi i \xi j}
\end{eqnarray*}

Thus $$\sum\limits_m f_a(x+m)e^{-2\pi i \xi m}\chi_{I_j}(x) = \hat{a}({\xi}) e^{2\pi i \xi j}  $$

Similarly $$\sum\limits_n g_b(x+n)e^{-2\pi i \eta  n}\chi_{I_j}(x) = \hat{b}({\eta}) e^{2\pi i \eta j}$$

Substituting these we get,
\begin{eqnarray*}
    \CC_\psi(f_a,g_b)(x)\chi_{I_j}(x) &=& \int_0^1\int_0^1 \psi(\xi,\eta)  \hat{a}(\xi)  \hat{b}(\eta)
    e^{2\pi i (\xi+\eta) j} d\xi d\eta \\
&=& \DD_\psi(a,b)(j),
\end{eqnarray*}
where $\DD_\psi$ is the bilinear operator
defined on  $l_{p_1}(\Z) \times l_{p_2}(\Z)$. Now
\begin{eqnarray*}
    |\DD_\psi(a,b)(j)|^{p_3} &=& |\CC_\psi(f_a,g_b)(x)|^{p_3}\chi_{I_j}, \\
 &=& 2 \int_{j+\frac{1}{4}}^{j+\frac{3}{4}} |\CC_\psi(f_a,g_b)(x)|^{p_3} dx
\end{eqnarray*}

\medskip
Summing over $j \in \Z$,
\begin{eqnarray*}
    \sum\limits_j |\DD_\psi(a,b)(j)|^{p_3} &=& 2 \sum\limits_j \int_{j+\frac{1}{4}}^{j+\frac{3}{4}} |\CC_\psi(f_a,g_b)(x)|^{p_3} dx \\
  &\leq 2&   \int_{\R}|\CC_\psi(f_a,g_b)(x)|^{p_3} dx \\
&\leq 2&   ||\CC_\psi||^{p_3}  ||f_a||_{p_1}^{p_3} ||g_b||_{p_2}^{p_3}
\end{eqnarray*}
Therefore
 $$||\DD_\psi(a,b)||_{p_3} \leq 2^{1/p_3}\ ||\CC_\psi|| \ ||a||_{p_1} \ ||b||_{p_2} .$$

\medskip

For the converse:-

\medskip

Let $\psi\in M_{{p_1},{p_2}}^{p_3}(\Z)$. Define
$$K_{n,m}:= \int_I \int_I \psi(\xi,\eta)e^{2 \pi i (\xi
n+\eta m)} d\xi d\eta$$ where $I = [0,1]$.

Let $f,g \in C_c^\infty(\R)$,we have ,
\begin{eqnarray*}
\CC_\psi(f,g)(x) &=& \int_{\R}\int_{\R}\hat{f}(\xi) \hat{g}(\eta)\psi(\xi,\eta)e^{2 \pi i x(\xi+\eta)} d\xi
d\eta\\
&=& \sum\limits_{n,m}
\int_I\int_I\hat{f}(\xi+n) \hat{g}(\eta+m)\psi(\xi+n,\eta+m)e^{2 \pi
i x(\xi+n+\eta+m)} d\xi d\eta\\
&=&\int_I \int_I \sum\limits_n \hat{f}(\xi+n)  e^{2 \pi i x(\xi+n)}\sum\limits_m \hat{g}(\eta+m)
e^{2 \pi i x(\eta+m)}\psi(\xi,\eta)e^{2 \pi i x(\xi+\eta)} d\xi d\eta\\
&=&\int_I\int_I \left( \sum\limits_n f(x+n)  e^{-2 \pi i
\xi n}\right) \left(\sum\limits_m g(x+m)e^{-2 \pi i \eta
m}\right) \psi(\xi,\eta) d\xi d\eta\\
&  = & \sum\limits_n \sum\limits_m f(x+n) g(x+m)
\int_I\int_I\psi(\xi,\eta)e^{-2 \pi i (\xi n+\eta m)}d\xi
d\eta\\
& =& \sum\limits_n \sum\limits_m K_{n,m}f(x-n) g(x-m)
\end{eqnarray*}

Hence,
$$ \int_{\R}|\CC_\psi(f,g)(x)|^{p_3} dx = \int_{\R}|\sum\limits_n \sum\limits_m K_{n,m}f(x-n) g(x-m)|^{p_3} dx $$
$$ = \int_I\sum\limits_l |\sum\limits_n \sum\limits_m K_{n,m}f(x+l-n) g(x+l-m)|^{p_3} dx$$
$$ \leq  \int_I ||\DD_\psi||^{p_3} \left(\sum\limits_n |f(x-n)|^{p_1} \right)^{\frac{p_3}{p_1}} \left(\sum\limits_m g(x-m)|^{p_2} \right)^{\frac{p_3}{p_2}} dx $$
Using the H\"{o}lder's inequality with the exponents
$\frac{p_1}{p_3},\frac{p_2}{p_3}$, we obtain
$$||\CC_\psi(f,g)||_{p_3} \leq ||\DD_\psi|| \ ||f||_{p_1} ||g||_{p_2}.$$
\qed

\medskip

We now turn our attention to periodic extensions of compactly supported bilinear multipliers.
We will prove the following result.
\medskip
%\subsection{\bf Proof of Theorem (\ref{compact support})}:
\begin{theorem}\label{compact support}

Let $\psi \in M_{{p_1},{p_2}}^{p_3}(\R)$ be such that $supp (\psi)
\subseteq I \times I$. where $I=[-1/2 ,1/2]$. Consider $\psi^\sharp$
the periodic extension of $\psi$ given by

$$\psi^\sharp (\xi,\eta) = \sum\limits_n \sum\limits_m \psi(\xi-n,\eta-m)$$
Then $\psi^\sharp \in M_{{p_1},{p_2}}^{p_3}(\R)$. Moreover, $||\psi^\sharp|| \leq c ||\psi||$,
where $c$ is a constant independent of $\psi$.
\end{theorem}
As a consequence of theorem (\ref{periodic symbols}) it would suffice to prove the following.
\begin{proposition}\label{1}

Let $\psi \in M_{{p_1},{p_2}}^{p_3}(\R)$ be such that $supp(\psi) \subseteq I \times I$.
where $I=[-1/2 ,1/2]$. Then $\psi^\sharp \in M_{{p_1},{p_2}}^{p_3}(\Z)$.
Moreover, $||\psi^\sharp|| \leq c ||\psi||$.
where $c$ is a constant independent of $\psi$.

\end{proposition}

We will need  the following two lemmas.
\begin{lemma}
Let $\psi \in M_{{p_1},{p_2}}^{p_3}(\R)$ be such that $supp(\psi) \subseteq I \times I$.
Then for $f,g \in \mathcal{S}(\R)$, $supp \ \widehat{\CC_\psi(f,g)} \subset [-2,2]$.
\end{lemma}
\noindent
{\bf Proof:} Let $h \in \mathcal{S}(\R)$ be such that $supp \ \hat{h}\subset[-2,2]^c$.
Then
\begin{eqnarray*}
 \left<\widehat{\CC_\psi(f,g)},h \right> &=& \left<\CC_\psi(f,g),\hat{h} \right>\\
&=& \int_{\R}\int_{\R}\int_{\R} \hat{f}(\xi) \hat{g}(\eta)\psi(\xi,\eta)e^{2 \pi i x(\xi+\eta)} d\xi
d\eta \hat{h}(x)dx \\
&=& \int_I\int_I \hat{f}(\xi) \hat{g}(\eta)\psi(\xi,\eta)\int_{\R}\hat{h}(x)e^{2 \pi i x(\xi+\eta)}dx d\xi d\eta \\
&=& \int_I\int_I \hat{f}(\xi) \hat{g}(\eta)\psi(\xi,\eta) h(\xi+\eta)d\xi d\eta  = 0
\end{eqnarray*}

Thus
$$\left<\widehat{\CC_\psi(f,g)},h \right> = 0$$
This proves the lemma.\qed

\medskip
For the proof of Proposition (\ref{1}) we will use the following result.

\begin{lemma}\label{AC}\cite{boas}
 Let $0 < p \leq \infty$ and $g$ be a slowly increasing $C^{\infty} $ function such that
 $supp(\hat{g})\subset [-R,R]$, then there exists a constant
$C> 0$, depending on $p$ such that
$$\sum\limits_n |g(n)|^p \leq C^p \  max(1,R)\int_{\R} |g(x)|^p dx $$
\end{lemma}

This is a well known sampling lemma.

\medskip

Now we prove Proposition (\ref{1}).

\medskip
\noindent
{\bf Proof:} Let $\Phi \in \mathcal{S}(\R)$.

Let $a = \{a_k\} $ and $b = \{b_l\} $ be two sequences with finitely many non-zero terms.

We define $f_a(x):= \sum\limits_k a_k \Phi(x-k) $ and $g_b(x):= \sum\limits_l b_l \Phi(x-l) $. It is easy to see that
$$||f_a||_{L^{p_1}(\R)} \leq c ||a||_{l^{p_1}(\Z)}$$
$$||g_b||_{L^{p_2}(\R)} \leq c ||b||_{l^{p_2}(\Z)}$$

where $c = \sum\limits_l |\Phi(x-l)| $.

Also, $\hat{f_a}(\xi) = \hat{\Phi}(\xi) \hat{a}(\xi) $. Similarly, \ $\hat{g_b}(\eta)= \hat{\Phi}(\eta) \hat{b}(\eta)$

We write the operator
\begin{eqnarray*}
\CC_\psi(f_a,g_b)(x)&=&\int_{\R}\int_{\R} \hat{f_a}(\xi) \hat{g_b}(\eta)\psi(\xi,\eta)e^{2 \pi i x(\xi+\eta)}
d\xi d\eta  \\
&=&\int_{I}\int_{I} \hat{\Phi}(\xi) \hat{a}(\xi) \hat{\Phi}(\eta) \hat{b}(\eta)\psi(\xi,\eta)
e^{2 \pi i x(\xi+\eta)} d\xi d\eta
\end{eqnarray*}

Choose $\Phi$ such that $\hat{\Phi}\equiv 1\ \ on \ \ I $. Then at the integer points we get

$$\CC_\psi(f_a,g_b)(n) = \int_{I}\int_{I}  \hat{a}(\xi)  \hat{b}(\eta)\psi(\xi,\eta)e^{2 \pi i n(\xi+\eta)}
d\xi d\eta =\DD_\psi(a,b)(n) $$

Using lemma (\ref{AC}) we get,
\begin{eqnarray*}
\sum\limits_n |\DD_\psi(a,b)(n)|^{p_3} = \sum\limits_n |\CC_\psi(f_a,g_b)(n)|^{p_3}& \leq& C_{p_3}^{p_3} 2^{p_3} \int_{\R}|\CC_\psi(f_a,g_b)(x)|^{p_3} dx\\
&\leq &  C_{p_3}^{p_3} 2^{p_3} ||\CC_\psi||^{p_3} ||f_a||_{p_1}^{p_3} ||g_b||_{p_2}^{p_3}\\
\end{eqnarray*}
i.e. $$||\DD_\psi(a,b)||_{p_3} \leq  C'_{p_3}||\CC_\psi|| \ ||a||_{p_1} ||b||_{p_2} $$
\qed
%\end{proof}
\section {\bf Jodeit type extension theorems}\label{jext}
In this section we will explore some extensions of bilinear multipliers
on $\T$  to bilinear multipliers on $\R$. Essentially our results
are analogues of Jodeit type of extensions in the linear case. Our
proofs are refinements of Jodeit's original proofs. For the sake of
completeness we include the proofs here. We will need the following
lemmas which may be of independent interest.

\medskip
In what follows $J$ will denote the interval $[-1/2, 1/2)$.

\begin{lemma}\label{p3} Let $\phi \in M_{p_1,p_2}^{p_3}(\T)$.
\begin{enumerate}
\item[(i)] If $p_3\geq1$ and  $a \in l^1(\Z^2)$  then
$a*\phi \in M_{p_1,p_2}^{p_3}(\T)$ and $||a*\phi|| \leq ||a||_1
||\phi||$.
\item[(ii)]If $p_3<1$ and  $a \in l^{p_3}(\Z^2)$  then
$a*\phi \in M_{p_1,p_2}^{p_3}(\T)$ and $||a*\phi|| \leq ||a||_{p_3}
||\phi||$.
\end{enumerate}
\end{lemma}
\noindent
{\bf Proof:} For $p_3 \geq 1$, this is an immediate consequence of Minkowski's inequality.
Assume $p_3< 1$ and Let $T'$ be the operator corresponding to $a*\phi$. 
For $f\in L^{p_1}(\T)$ and $g\in L^{p_2}(\T)$,
$$||T'(f,g)||_{p_3}^{p_3} = \int_J |\sum\limits_{n,m} a* \phi(n,m) \hat{f}(n) \hat{g}(m) \ e^{2 \pi i x (n+m)}|^{p_3} dx$$
$$= \int_{J} |\sum\limits_{l,k}\sum\limits_{n,m} a(l,k)\phi(n-l,m-k) \hat{f}(n) \hat{g}(m) \ e^{2 \pi i x (n+m)}|^{p_3} dx$$
$$\leq \int_{J}\sum\limits_{l,k}|a(l,k)|^{p_3}\ \ |\sum\limits_{n,m}\phi(n-l,m-k) \hat{f}(n) \hat{g}(m) \ e^{2 \pi i x (n+m)}|^{p_3} dx$$
$$\leq \|T\|^{p_3} \|a\|_{p_3}^{p_3}  \|f\|_{p_1}^{p_3} \|g\|_{p_2}^{p_3}$$

The first inequality follows from $|\sum\limits_i \alpha_i|^p \leq \sum\limits_i|\alpha_i|^p$, $ 0<p<1$
and the second from the assumption that the operator $T$ is bounded. Hence we obtain
$$\|T'(f,g)\|_{p_3} \leq \|T\| \|a\|_{p_3}  \|f\|_{p_1} \|g\|_{p_2} $$\qed

\begin{lemma}\label{dilation} Let $\phi \in M_{p_1,p_2}^{p_3}(\T)$. For a positive integer $k$ ,we define $\phi_k$ as follows
$\phi_k(n,m):= \phi (n/k,m/k)$ if $k$ divides both $n,m $, and $\phi_k(n,m):=0$
otherwise. Then $\phi_k\in M_{p_1,p_2}^{p_3}(\T)$ with norm not exceeding that of $\phi$.
\end{lemma}\
\noindent {\bf Proof:} For $f\in L^\infty(\T)$, let $F(x)= \frac{1}{k}\sum\limits_0^{k-1}
f(\frac{x+j}{k})$.
\begin{eqnarray*}
\hat{F}(n)&=& \int_J \frac{1}{k}\sum\limits_0^{k-1} f(\frac{x+j}{k})e^{-2 \pi i x.n}dx\\
&=& \frac{1}{k}\sum\limits_0^{k-1}\int_J f(\frac{x+j}{k})e^{-2 \pi i x.n}dx\\
&=& \int_J f(x)e^{-2 \pi i x.kn}dx\\
&=& \hat{f}(kn) \end{eqnarray*}
Also note that
$\|F\|_1 \leq \|f\|_1$ and $\|F\|_{\infty} \leq \|f\|_{\infty}$.
Hence $$\|F\|_p \leq \|f\|_p,    \ \ 1 \leq p \leq \infty$$

Similarly for $g\in L^\infty(\T)$, we define $G$. Let $T_k$
be the operator corresponding to $\phi_k$. Then
\begin{eqnarray*}
T_k(f,g)(x) &=& \sum\limits_n \sum\limits_m \hat{f}(n)\hat{g}(m)\phi_k(n,m)
e^{2\pi i
 x(n+m)}\\
&=& \sum\limits_n \sum\limits_m \hat{f}(kn)\hat{g}(km)\phi(n,m) e^{2\pi i
 x(kn+km)}\\
&=& \sum\limits_n \sum\limits_m \hat{F}(n)\hat{G}(m)\phi(n,m) e^{2\pi i
 x(kn+km)}\\
 &=& T(F,G)(kx)
\end{eqnarray*}
Hence \begin{eqnarray*}
\|T_k(f,g)\|_{p_3}^{p_3} &=& \int_J |T(F,G)(kx)|^{p_3}dx\\
&=& \frac{1}{k} \int_{-k/2}^{k/2} |T(F,G)(x)|^{p_3}dx\\
&=& \int_J |T(F,G)(x)|^{p_3}dx\\
&\leq& \|T\|^{p_3} \|F\|_{p_1}^{p_3} \|G\|_{p_2}^{p_3}\\
&\leq& \|T\|^{p_3} \|f\|_{p_1}^{p_3} \|g\|_{p_2}^{p_3}
\end{eqnarray*}
Thus we obtain that $\phi_k$ is in $M_{p_1,p_2}^{p_3}(\T)$.\qed
\medskip

Our first extension result is the following theorem.
\begin{theorem}\label{extension1} Let  $\phi$ be in  $M_{p_1,p_2}^{p_3}(\T)$ and  $S$ be a function  supported in $\frac{1}{2} J \times \frac{1}{2} J $ such that  its periodic extension $S^{\sharp}$ from $J \times J$ satisfies $\sum\limits_{n,m}|\hat{S^{\sharp}}(n,m)|^p< \infty $ , where $p =$ min $(1, p_3)$.\\
Then $\psi(\xi,\eta):= \sum\limits_{n,m} \phi (n,m) \hat{S}(\xi-n,\eta-m)\in M_{p_1,p_2}^{p_3}(\R)$. Moreover, $\|\psi\| \leq c \|\phi\|$, where $c = 2^{1/p}\sum\limits_{n,m} |\hat{S^{\sharp}}(n,m)|^p< \infty $.
\end{theorem}
\noindent {\bf Proof:} It is enough to prove that
$\psi_{r}(\xi,\eta) = \sum\limits_{l,k} \phi (l,k)\  r^{|l|+|k|}
\hat{S}(\xi-l,\eta-k)$ belongs to $M_{p_1,p_2}^{p_3}(\R)$ for
$0<r<1$  with $\|C_{\psi_r}\|\leq c\|\phi\|.$ Let $k_{\phi_r}$ be
the kernel corresponding to the bi-linear multiplier
$\phi(n,m)r^{|n|+|m|}$. Clearly, $\widehat{K_{\phi_r}S}=\psi_r$, considered as a function on $\R^2$. From lemma (\ref{p3})
$\widehat{K_{\phi_r}S^\#}=\hat K_{\phi_r}* \hat{S^\#}$ belongs to
$M_{p_1,p_2}^{p_3}(\T).$

 Let $f\in{\mathcal S}(\R)$ and  for each $n\in\Z$, let $f_n$ denote the  1-periodic extension of the function
 $f(x+n/2)\chi_J (x)$ from $J$.
Then it can be easily verified that
\begin{eqnarray*}\sum\limits_n\|f_n\|_{L^p(\T)}^p &\leq& 2\|f\|_p^p
\end{eqnarray*}
Now  for $x\in\frac{1}{2}J$ we have,
\begin{eqnarray*}
C_{\psi_r}(f,g)(x+\frac{n}{2})&=&\int_{\frac{1}{2}J}\int_{\frac{1}{2}J}f(x-t+\frac{n}{2})g(x-s+\frac{n}{2})(K_{\phi_r}S)(t,s)
dt ds\\
&=& \int_{J}\int_{J}f_n(x-t)g_n(x-s)(K_{\phi_r}S^\#)(t,s) dt ds
\end{eqnarray*}
Thus,
\begin{eqnarray*}
\|C_{\psi_r}(f,g)\|_{p_3}^{p_3} & \leq &  \sum\limits_n
(c\|\phi\|)^{p_3}\|f_n\|_{p_1}^{p_3}\|g_n\|_{p_2}^{p_3}\\
&\leq & 2(c\|\phi\|)^{p_3}\|f\|_{p_1}^{p_3}\|g\|_{p_2}^{p_3}
\end{eqnarray*}
The last inequality follows as an application of H\"{o}lder's
inequality.
\qed
\medskip

Our next result is the piecewise linear extension of $\phi\in M_{p_1,p_2}^{p_3}(\T)$ to $\psi\in M_{p_1,p_2}^{p_3}(\R)$.

\begin{theorem}\label{thm} Let $\phi \in M_{p_1,p_2}^{p_3}(\T)$ and $p_3>\frac{1}{2}$. If
$\Lambda(x_1,x_2) = (1-|x_1|)(1-|x_2|)$ in $[-1,1)\times [-1,1)$
and\ \ $0$ otherwise. Then the function $\psi$ defined as
$\psi(\xi,\eta):= \sum\limits_{n,m} \phi (n,m) \Lambda(\xi-n,\eta-m) \in
M_{p_1,p_2}^{p_3}(\R)$ and $||\psi|| \leq c  ||\phi|| $.
\end{theorem}
\noindent {\bf Proof:} Consider $S(x,y)=\frac {sin^2 4\pi x}{(4\pi
x)^2}\frac {sin^2 4\pi y}{(4\pi y)^2}$. Let $S_{k,l}(x,y)$ be the 1 - periodic extension of $\chi_{\frac{k}{2}+{\frac{1}{2}J}}(x)\chi_{\frac{l}{2}+{\frac{1}{2}J}}(y) S(x,y)$ from $(\frac{k}{2}+J)\times (\frac{l}{2}+J)$.
An easy computation using integration by parts gives that for $n,m \neq 0$ $$|\hat
S_{k,l}(n,m)|\leq \frac{C}{(1+k^2)(1+n^2)(1+l^2)(1+m^2)}.$$ Hence,
by Theorem \ref{extension1} we have $\sum\limits_{n,m} \phi (n,m) \hat
S_{k,l}(\xi-n,\eta-m)$ belongs to $M_{p_1,p_2}^{p_3}(\R)$ with
bounds not exceeding $\frac{C}{(1+k^2)(1+l^2)}\|\phi\|$.\\
Thus $\sum\limits_{n,m} \phi (n,m) \hat
S(\xi-n,\eta-m)= \sum\limits_{n,m} \phi (n,m)
\Lambda(\frac{\xi-n}{4},\frac{\eta-m)}{4})$ is in
$M_{p_1,p_2}^{p_3}(\R)$. By applying Lemma \ref{dilation} we get
$\psi\in M_{p_1,p_2}^{p_3}(\R)$ with the required bound.\qed

\medskip

As a consequence of this we get the desired piecewise constant extension result i.e.,

\begin{theorem}\label{piecewise}Let $\phi$ be in $M_{p_1,p_2}^{p_3}(\T)$, where $p_1, p_2 > 1$.
Then $\psi(\xi,\eta):= \sum\limits_{n,m} \phi (n,m) \chi_{J \times
J}(\xi-n,\eta-m) \in M_{p_1,p_2}^{p_3}(\R)$.
\end{theorem}
\noindent {\bf Proof:} Define  $\phi_2(n,m):=\phi(n/2, m/2)$
if $n,m$ both are even and $\phi_2(n,m):=0$ otherwise. By Lemma (\ref{dilation}),
$\phi_2(n,m) \in M_{p_1,p_2}^{p_3}(\T)$ and $||\phi_2|| \leq
||\phi||$. Consider $\theta_2(n,m) = \phi_2(n,m) + \phi_2(n-1,m)+\phi_2(n,m-1) + \phi_2(n-1,m-1)$.
Clearly, $\theta_2 \in M_{p_1,p_2}^{p_3}(\T)$ and $||\theta_2|| \leq
4 ||\phi||$.

Let $\Lambda(\xi, \eta) = (1-|\xi|)(1-|\eta|)$ in $[-1,1)\times
[-1,1)$, and $0$ elsewhere. Then by Theorem \ref{thm},
$\Theta_2(\xi,\eta)= \sum\limits_{n,m} \theta_2(n,m)\Lambda(\xi-n,\eta-m)$
is in $ M_{p_1,p_2}^{p_3}(\R)$. Note that $\Theta_2(\xi,\eta)=
\phi(n,m)$ if $(\xi, \eta)\in [2n,2n+1]\times [2m, 2m+1].$ Let
$\tilde{\chi}$ be the 2-periodic extension of a function which is $1$ for
$0<x<1$ and $0$ for $-1<x<0$. We know that $ \tilde{\chi} \in M_p(\R)$ for
$p>1$. Hence $\Psi_2(\xi, \eta) = \tilde{\chi}(\xi)\tilde{\chi} (\eta) \Theta_2(\xi,
\eta) \in M_{p_1,p_2}^{p_3}(\R)$. Now $\Psi_2(\xi, \eta) =
\phi(n,m)$ for $2n<\xi<2n+1$, $2m<\eta<2m+1$ and $\Psi_2(\xi,
\eta)=0$ otherwise. Since $\psi(\xi, \eta)= \Psi_2(2\xi, 2\eta)+\Psi_2(2\xi+1,
2\eta)+\Psi_2(2\xi, 2\eta+1)+\Psi_2(2\xi+1,2\eta+1)$. The result
follows. \qed

\medskip
{\bf Remark}: The above theorem does not hold if either of $p_1,
p_2$ is $1$. This is very easy to verify. Without loss of generality, we can assume that $p_1 = 1$.
Let $\tilde{\phi} \in M_1(\T)$, and $T \in M_1(\R)$  corresponding to the piece-wise constant extension of $\tilde{\phi}$. Put $\phi (n,m)= \tilde{\phi} (n)$. Then by H\"{o}lder's
inequality $\phi\in M_{1,p_2}^{p_3}(\R)$. Suppose the piecewise constant extension $\psi(\xi,\eta)= \sum\limits_{n,m} \phi (n,m) \chi_{J \times J}(\xi-n,\eta-m)= \sum\limits_n \tilde{\phi} (n) \chi_{J}(\xi-n)$ is in $M_{1,p_2}^{p_3}(\R)$. Then for  $f \in L_1(\R)$ and $g\in L_{p_2}(\R)$, we have $\|C_\psi(f,g)\|_{p_3}= \|T(f).g \|_{p_3}\leq c \|f\|_1 \|g\|_{p_2}$. Now notice that $\frac{1}{p_3}$ and $\frac{p_2}{p_3}$ are conjugate indices and $\|f\|_1^{p_3} = \| |f|^{p_3} \|_{\frac{1}{p_3}}$. Hence by using duality we get $\|Tf\|_{1} \leq c \|f\|_{1}$. But this is not true for any nonconstant $\tilde{\phi}$.


\begin{thebibliography}{1}

\bibitem{auscher} Auscher, P.; Carro, M. J., {\it On relations between operators on $R\sp N,\;T\sp N$ and $Z\sp N$,} Studia Math. 101 (1992), no. 2, 165--182.
\bibitem{blasco} Blasco, O., {\it Bilinear multipliers and transference,} I.J.M.M.S. 2005, No 4 (2005), 545-554.

\bibitem{bcg} Blasco, O. , Carro, M. , Gillespie, T. A., {\it Bilinear Hilbert Transform on measure spaces,} J. Fourier Anal. appl. 11(2005), no. 4, 459--470.

\bibitem{boas}Boas, R.P., {\it Entire Functions,} Academic Press, 1954.


\bibitem{dlw} de Leeuw, K.,
{\it On $L_p$ multipliers,} Ann. of Math. 81 (1965), 364--379.
\bibitem{fs} Fan, Dashan, Sato, Shuichi, {\it  Transference on certain  multilinear multiplier operators,} J. Aust.  Math. Soc. 70 (2001), no. 1,  37--55.

\bibitem{gtor} Grafakos, L., Torres, R. H., {\it Multilinear Calder$\acute{o}$n-Zygmund theory,} Advances in Mathematics 165 (2002), no. 1, 124--164.

\bibitem{ghoz} Grafakos, L., Honzik, P., {\it Maximal transference and summability of multilinear Fourier series,} J. Aust. Math. Soc., 80 (2006), no. 1, 65--80.

\bibitem{jod} Jodeit, M.A.,
{\it Restrictions and extensions of Fourier multipliers,} Studia Math. 34 (1970) , 215--226.

\bibitem{lt1} Lacey, Michael, Thiele, Christoph,
{\it $L^p$ estimates on the bilinear Hilbert transform for $2<p<\infty$,} Ann. of Math.(2) 146 (1997), no. 3, 693--724.

\bibitem{lt2} Lacey, Michael, Thiele, Christoph,
{\it On Calder$\acute{o}$n's conjecture,} Ann. of Math.(2) 149 (1999), no. 2, 475--496.

\bibitem{sp} Madan, Shobha; Mohanty, Parasar, {\it Jodeit's extension for bilinear multipliers,} {\it to appear in Bull. London Math Soc.}
\end{thebibliography}
\end{document}